\documentclass{aip-cp}
\usepackage[numbers]{natbib}
\usepackage{graphicx}
\usepackage{color, multirow,amsfonts}
\usepackage{url}
\usepackage{rotating}
\begin{document}
\title{Logarithm of ratios of two order statistics\\ and regularly varying tails}
\author[aff1]{Pavlina K. Jordanova\corref{cor1}}
\author[aff3,aff4]{Milan Stehl\'\i k}

\affil[aff1]{Faculty of Mathematics and Informatics, Konstantin Preslavsky University of Shumen, \\115 "Universitetska" str., 9712 Shumen, Bulgaria.}
\affil[aff3]{Institute of Statistics, Universidad de Valpara\'iso, Valpara\'iso, Chile.}
\affil[aff4]{Department of Applied Statistics and Linz Institute of Technology, Johannes Kepler University, Altenbergerstrasse 69, 4040 Linz, Austria.}

\corresp[cor1]{Corresponding author: pavlina\_kj@abv.bg}

\maketitle

\begin{abstract} Here we suppose that the observed random variable has cumulative distribution function $F$ with regularly varying tail, i.e. $1-F \in RV_{-\alpha}$, $\alpha > 0$. Using the results about exponential order statistics we investigate logarithms of ratios of two order statistics of a sample of independent observations on Pareto distributed random variable with parameter $\alpha$. Short explicit formulae for its mean and variance are obtained. Then we transform this function in such a way that to obtain unbiased, asymptotically efficient, and asymptotically normal estimator for $\alpha$. Finally we simulate Pareto samples and show that in the considered cases the proposed estimator outperforms the well known Hill, t-Hill, Pickands and Deckers-Einmahl-de Haan estimators.
\end{abstract}

\section{HISTORY OF THE PROBLEM}

The usefulness of regularly varying (RV) functions in economics seems to be discussed for the first time during modeling of the wealth in our society by Pareto distribution, called to the name of Vilfredo Pareto (1897). J. Karamata (1933) provides their definition and integral representation. Later on the Convergence to types theorem, proved by R. A. Fisher, L. H. C. Tippett (1928), and B.V. Gnedenko (1948) plays a key role for their future applications. It is well known that this class of distributions describes very well the domain of attraction of stable distribution (see Mandelbrot (1960) \cite{Mandelbrot1960}) and max-domain of attraction of Fr$\acute{e}$chet distribution (see M. Fr$\acute{e}$chet (1927)). Laurens de Haan (1970) and co-authors \cite{dH70,deHaanStadtmueller,deHaanFerreira} develop the main machinery for working with  cumulative distribution functions(c.d.fs.) with such tail behaviour. Let us remind that the c.d.f. $F$ has regularly varying right tail with parameter $\alpha > 0$, if
$$\lim_{t\to \infty} \frac{1-F(xt)}{1-F(t)}= x^{-\alpha}, \quad \forall x > 0.$$
After their works the topic spread over the world very fast and many estimators of the index of regular variation are proposed, see e.g. Hill (1975) \cite{Hill}, Pickands (1975)\cite{Pickands} and Deckers-Einmahl-de Haan (1989) \cite{Dekkers1989}, t-Hill (Stehlik and co-authors (2010) \cite{Stehlik2010,Fabian,Stehlik2012}, and Pancheva and Jordanova (2012) \cite{JordanovaPancheva, JordanovaMilan2012}), among others.

Here we show the usefulness of functions of two central order statistics in estimating the parameter of regular variation. Under very general settings we show that the logarithm of the fraction of two specific central order statistics is an weakly consistent and asymptotically normal estimators of the logarithm of the corresponding theoretical quantiles. Then we use these functions and obtain our estimator for $\alpha$. Its main advantage is that it is very flexible and provides a useful accuracy given mid-range and small samples.
Pareto case, considered in Section 3 motivates our investigation. First we define a biased form of the estimator. Then using results about order statistics, which could be seen e.g. in Nevzorov (2001) \cite{Nevzorov} we obtain explicit formulae for its mean and variance. This allows us to define unbiased correction which is asymptotically efficient. Then we prove asymptotic normality and obtain large sample confidence intervals. Our simulation study depicts the advantages of the considered estimators over Hill, t-Hill, and Deckers-Einmahl-de Haan estimators.  The paper finishes with some conclusive remarks.

Trough the paper we assume that $\mathbf{X}_1, \mathbf{X}_2, ..., \mathbf{X}_n$ are independent observations on a random variable(r.v.) $\mathbf{X}$, and denote by $\mathbf{X}_{(1, n)} \leq \mathbf{X}_{(2, n)} \leq ... \leq \mathbf{X}_{(n, n)}$ the corresponding increasing order statistics.
$$H_{n, m} =  1 + \frac{1}{2^m} + \frac{1}{3^m} + ... + \frac{1}{(n-1)^m} + \frac{1}{n^m}, \quad n = 1, 2, ..., $$
denotes the $n$-the Generalized harmonic number of power $m = 1, 2, ...$, and $H_n = H_{n, 1}$, $n = 1, 2, ...$ is for the well-known $n$-th harmonic number.

The main object of interest in this point are the statistics
$$ Q_{k,s} : = \frac{log\frac{\mathbf{X}_{(ks, (s+1)k-1)}}{\mathbf{X}_{(k, (s+1)k-1)}}}{H_{ks-1} - H_{k-1}}, \quad Q_{k,s}^* = \frac{log\frac{\mathbf{X}_{(ks, (s+1)k-1)}}{\mathbf{X}_{(k, (s+1)k-1)}}}{log(s)}, \quad s = 2, 3, ...$$
The estimator $Q_{k,3}^*$ it is obtained in Jordanova et al. \cite{jordanova2017measuring} via quantile matching procedure. About the last procedure see e.g. Sgouropoulos et al. (2015) \cite{sgouropoulos2015matching}.

Along the paper $\stackrel{d}{\to}$ means convergence in distribution.

\section{GENERAL RESULTS}

In 1933 - 1949 Smirnoff \cite{smirnov1949limit} shows that in case of central order statistics, and more precisely for $k, n$ and $p$ such that $\frac{k}{n} \to p \in (0, 1)$ and  $\sqrt{n}\left(\frac{k}{n} - p\right) \to \mu \in (-\infty, \infty)$, the asymptotic distribution of $\frac{\sqrt{n}[X_{(k, n)} - F^\leftarrow(p)]}{\sqrt{\frac{p(1-p)}{f^2[F^\leftarrow(p)]}}}$ is a standard normal. Moreover it seems that he has a similar results about bivariate order statistics. It could be seen e.g. in Arnold et al. (1992) \cite{arnold1992first}, p. 226, Mosteller (1946) \cite{Mosteller1946} p.338, Nair \cite{Nair2013}, p.330, or Wilks \cite{Wilks1948} among others. The multivariate delta method is a very powerful technique for obtaining confidence intervals in such cases. In the next theorem we apply them and obtain the limiting distribution of the logarithmic differences of central order statistics.

{\bf Smirnoff's theorem.} Assume for $n \to \infty$, $\frac{k_1(n)}{n} \to p_1 \in (0, 1)$, $\frac{k_2(n)}{n} \to p_2 \in (0, 1)$, $f\left[F^\leftarrow(p_1)\right] \in (0, \infty)$,  and $f\left[F^\leftarrow\left(p_2\right)\right] \in (0, \infty)$. Then
$$\sqrt{n}\left(
    \begin{array}{c}
      X_{(k_1, n)} - F^\leftarrow(p_1) \\
      X_{(k_2, n)} - F^\leftarrow(p_2)\\
    \end{array}
  \right)
\stackrel{d}{\to} \left(
                    \begin{array}{c}
                      \theta_1 \\
                      \theta_2 \\
                    \end{array}
                  \right), \quad \left(
                    \begin{array}{c}
                      \theta_1 \\
                      \theta_2 \\
                    \end{array}
                  \right) \in N(0, V)
$$
where the covariance matrix
$$V = \left(
        \begin{array}{cc}
          \frac{p_1(1-p_1)}{f^2[F^\leftarrow(p_1)]} & \frac{p_1(1-p_2)}{f^2[F^\leftarrow(p_1)]f^2[F^\leftarrow(p_2)]} \\
         \frac{p_1(1-p_2)}{f^2[F^\leftarrow(p_1)]f^2[F^\leftarrow(p_2)]} & \frac{p_2(1-p_2)}{f^2[F^\leftarrow(p_2)]} \\
        \end{array}
      \right).
$$

We apply this theorem together with the Multivariate delta method and obtain asymptotic normality of the estimators, discussed in this paper.

{\bf Theorem 1.} Consider a sample of $n = (s+1)k-1$, $s = 2, 3, ...$ independent observations on a r.v. $X$ with c.d.f. $F$ and p.d.f. $f = F'$. If there exists $0 < f\left[F^\leftarrow(\frac{1}{s+1})\right] < \infty$ and $0 < f\left[F^\leftarrow(\frac{s}{s+1})\right] < \infty$, then for $k \to \infty$
\begin{equation}\label{GeneralAN}
T_{k,s} := \sqrt{(s+1)k-1}\left[log\left(\frac{X_{(ks, (s+1)k-1)}}{X_{(k, (s+1)k-1)}}\right) -  log\left(\frac{F^\leftarrow(\frac{s}{s+1})}{F^\leftarrow(\frac{1}{s+1})}\right)\right] \stackrel{d}{\to } N\left( 0; V\right).
 \end{equation}
The variance $V$ in (\ref{GeneralAN}) is $V = \frac{1}{(s+1)^2} \left(\frac{s}{a_{F,s}^2} - \frac{2}{a_{F,s}b_{F,s}} + \frac{s}{b_{F,s}^2}\right)$, where
$a_{F,s} = F^\leftarrow(\frac{1}{s+1})f\left[F^\leftarrow(\frac{1}{s+1})\right] = \frac{1}{\{log[ F^\leftarrow(p)]\}'|_{p = \frac{1}{s+1}}}$, and $b_{F,s} = F^\leftarrow(\frac{s}{s+1})f\left[F^\leftarrow(\frac{s}{s+1})\right] =$ $ \frac{1}{\{log[ F^\leftarrow(p)]\}'|_{p = \frac{s}{s+1}}}.$

{\bf Proof:} We will apply the Theorem of Smirnoff for $p_1 = \frac{1}{s+1}$ and $p_2 = \frac{s}{s+1}$ and Multivariate delta method.

By assumptions the conditions $f\left[F^\leftarrow\left(\frac{1}{s+1}\right)\right] \in (0, \infty)$, $f\left[F^\leftarrow\left(\frac{s}{s+1}\right)\right] \in (0, \infty)$ are satisfied. And for $k \to \infty$ we have $\frac{sk}{(s+1)k-1} \to \frac{s}{s+1}$, $\frac{s}{(s+1)k-1} \to \frac{1}{s+1}$, $\sqrt{(s+1)k-1}\left(\frac{sk}{(s+1)k-1} - \frac{s}{s+1}\right) \to 0$ and $\sqrt{(s+1)k-1}\left(\frac{k}{(s+1)k-1} - \frac{s}{s+1}\right) \to 0$, therefore the Smirnoff's theorem on the joint asymptotic normality of the order statistics, says that
 $$\sqrt{(s+1)k-1}\left(\begin{array}{c}
                             X_{(k, (s+1)k-1)} - F^\leftarrow(\frac{1}{s+1}) \\
                             X_{(ks, (s+1)k-1)} - F^\leftarrow(\frac{i}{s+1}) \\
                           \end{array}
                         \right) \stackrel{d}{\to} N\left[\left(\begin{array}{c}
                         0\\
                         0 \\
                           \end{array}\right); D\right], \quad k \to \infty,$$
 where the asymptotic covariance matrix of this bivariate distribution is
 $$D = \frac{1}{(s+1)^2}\left(\begin{array}{cc}
 \frac{s}{f^2\left[F^\leftarrow(\frac{1}{s+1})\right]} & \frac{1}{f\left[F^\leftarrow(\frac{1}{s+1})\right]f\left[F^\leftarrow(\frac{s}{s+1})\right]} \\
 \frac{1}{f\left[F^\leftarrow(\frac{1}{s+1})\right]f\left[F^\leftarrow(\frac{s}{s+1})\right]} & \frac{s}{f^2\left[F^\leftarrow(\frac{s}{s+1})\right]} \\
 \end{array}\right)$$
 and the asymptotic correlation between these two order statistics is $\frac{1}{s}$.

 Consider the function $g(x, y) = log\left(\frac{y}{x}\right)$. For $x > 0$ and $y > 0$ it is continuously differentiable.

The Jacobian of the transformation is
$$J : = \left[\frac{\partial g(x, y)}{\partial x}, \frac{\partial g(x, y)}{\partial y}\right] = \left(-\frac{1}{x}, \frac{1}{y}\right).$$

 The asymptotic mean is
 $$ \lim_{k \to \infty} E log\left(\frac{X_{(ks, (s+1)k-1)}}{X_{(k, (s+1)k-1)}}\right) = g\left[F^\leftarrow(\frac{1}{s+1}), F^\leftarrow(\frac{s}{s+1})\right] = log\left(\frac{F^\leftarrow(\frac{s}{s+1})}{F^\leftarrow(\frac{1}{s+1})}\right).$$

Now we apply the Multivariate Delta method, which could be seen e.g. in Sobel (1982) \cite{MultivariateDeltaMethod}, and obtain that the asymptotic variance of $T_{k,s}$ is
 \begin{eqnarray*}
  V : &=& J \times D \times J' = \left(-\frac{1}{x}, \frac{1}{y}\right)\left|_{x = F^\leftarrow(\frac{1}{s+1}), y = F^\leftarrow(\frac{s}{s+1})}\right. \\
    &\times&  \frac{1}{(s+1)^2}\left(\begin{array}{cc}
 \frac{s}{f^2\left[F^\leftarrow(\frac{1}{s+1})\right]} & \frac{1}{f\left[F^\leftarrow(\frac{1}{s+1})\right]f\left[F^\leftarrow(\frac{s}{s+1})\right]} \\
 \frac{1}{f\left[F^\leftarrow(\frac{1}{s+1})\right]f\left[F^\leftarrow(\frac{s}{s+1})\right]} & \frac{s}{f^2\left[F^\leftarrow(\frac{s}{s+1})\right]} \\
 \end{array}\right)  \left(
                             \begin{array}{c}
                              -\frac{1}{x} \\
                               \frac{1}{y} \\
                             \end{array}
                           \right)\left|_{x = F^\leftarrow(\frac{1}{s+1}), y = F^\leftarrow(\frac{s}{s+1})}\right. \\
   &=& \frac{1}{(s+1)^2} \left[\left(-\frac{1}{x}, \frac{1}{y}\right)\left(\begin{array}{cc}
 \frac{s}{f^2(x)} & \frac{1}{f(x)f(y)} \\
 \frac{1}{f(x)f(y)} & \frac{s}{f^2(y)} \\
 \end{array}\right) \left(
                             \begin{array}{c}
                              -\frac{1}{x} \\
                               \frac{1}{y} \\
                             \end{array}
                           \right)\right]\left|_{x = F^\leftarrow(\frac{1}{s+1}), y = F^\leftarrow(\frac{s}{s+1})}\right. \\
    &=& \frac{1}{(s+1)^2} \left[\frac{s}{x^2f^2(x)} - \frac{2}{xyf(x)f(y)} + \frac{s}{y^2f^2(y)}\right] \left|_{x = F^\leftarrow(\frac{1}{s+1}), y = F^\leftarrow(\frac{s}{s+1})}\right. =\frac{1}{(s+1)^2} \left(\frac{s}{a_{F,s}^2} - \frac{2}{a_{F,s}b_{F,s}} + \frac{s}{b_{F,s}^2}\right).
 \end{eqnarray*}
 \hfill Q.A.D.

Slutsky's theorem about continuous functions together with the definition of convergence in probability, application of quantile transform, and Smirnoff's theorem about a.s. convergence of empirical quantiles to corresponding theoretical one, lead us to the following result. Without lost of generality we consider only a.s. positive r.vs, however the result could be easily transformed for $P(X < \mu) = 1$ or $P(X > \mu) = 1$, $\mu > 0$.

{\bf Theorem 2.} Assume $P(X > 0) = 1$. If $f\left[F^\leftarrow\left(\frac{1}{s+1}\right)\right] \in (0, \infty)$, $f\left[F^\leftarrow\left(\frac{s}{s+1}\right)\right] \in (0, \infty)$, then for $s = 2, 3, ...$
\begin{equation}\label{logfraction}
 log\left(\frac{X_{(ks, (s+1)k-1)}}{X_{(k, (s+1)k-1)}}\right) \stackrel{P}{\to } log\left[\frac{F^\leftarrow(\frac{s}{s+1})}{F^\leftarrow(\frac{1}{s+1})}\right], \quad k \to \infty.
 \end{equation}

\section{PARETO CASE}

In this section we assume that $\mathbf{X}_1, \mathbf{X}_2, ..., \mathbf{X}_n$ are independent observations on a r.v. $\mathbf{X}$ with Pareto c.d.f.
\begin{equation}\label{Pareto}
 F_{\mathbf{X}}(x) = \left\{
                  \begin{array}{ccc}
                    0 & , & x \leq \delta \\
                    1 - \left(\frac{\delta}{x}\right)^{\alpha} & , & x > \delta
                  \end{array}
                \right., \quad \alpha > 0, \quad \delta > 0.
\end{equation}
Briefly we will denote this by $\mathbf{X} \in Par(\alpha, \delta)$.  Different generalizations of this distributions could be seen in Arnold (2015) \cite{Arnold2015}. The number $-\alpha$ is called "index of regular variation of the tail of c.d.f.". It determines the tail behaviour of the c.d.f. See  e.g. de Haan and Ferreira \cite{deHaanFerreira}, Resnick \cite{Resnick87}, or Jordanova \cite{Jordanova2019}.

Denote by $\mathbf{X} \in Exp(\lambda)$, $\lambda > 0$ the fact that the r.v. $\mathbf{X}$ has c.d.f.
\begin{equation}\label{Exponential}
 F_{\mathbf{X}}(x) = \left\{
                  \begin{array}{ccc}
                    0 & , & x \leq 0 \\
                    1 - e^{-\lambda x} & , & x > 0
                  \end{array}
                \right..
\end{equation}
The results in the following theorem allow us later on, in Corollaries 1 and 2, to obtain unbiased, consistent,  and asymptotically efficient estimators of the parameter $\alpha$.

{\bf Theorem 3.} Assume $\mathbf{X}_{(1, n)} \leq \mathbf{X}_{(2, n)} \leq ... \leq \mathbf{X}_{(n, n)}$ are order statistics of independent observations on a r.v. $\mathbf{X} \in Par(\alpha, \delta)$, $\alpha > 0$, $\delta > 0$, and $1 \leq i < j \leq n$ are integer.
 \begin{description}
 \item[i)] Denote by $\rho$ a Beta distributed with parameters $n - j + 1$, and $j - i$. Then
 $$log\left(\frac{X_{(j, n)}}{X_{(i, n)}}\right)  \stackrel{d}{=} - \frac{1}{\alpha} \log(\rho) \stackrel{d}{=} E_{(j - i, n-i)} \stackrel{d}{=} \frac{1}{\alpha} E_{(j - i, n-i)}^*,$$
where $E_{(j - i, n - i)}$ is the $j - i$-th order statistics in a sample of $n - i$ independent observations on i.i.d. Exponential r.vs. with parameter $\alpha$, and $E_{(j - i, n - i)}^*$ is the $j - i$ - th order statistic of a sample of $n-i$ independent observations on exponentially distributed r.v. with parameter $1$. Its probability density function is
 $$f_{log\frac{\mathbf{X}_{(j, n)}}{\mathbf{X}_{(i, n)}}}(x) = \frac{\alpha (n-i)!}{(j-i-1)!(n-j)!}(1-e^{-\alpha x})^{j-i-1}e^{-\alpha x(n-j + 1
 )}, \quad x > 0.$$
 \item[ii)] $\mathbb{E} \left[log\frac{\mathbf{X}_{(j, n)}}{\mathbf{X}_{(i, n)}}\right] = \frac{1}{\alpha}\left(H_{n-i} - H_{n-j}\right),$ and
   $\mathbb{D} \left[log\frac{\mathbf{X}_{(j, n)}}{\mathbf{X}_{(i, n)}}\right] = \frac{1}{\alpha^2}\left(H_{n-i, 2} - H_{n-j, 2}\right).$
 \end{description}

{\bf Proof:}  Let us fix  $1 \leq i < j \leq n$, integers.
Because of $g(x) = e^x$ is a strictly increasing function, it is well known that the probability quantile transform, entails
$$ \left(\frac{\mathbf{X}_{(1, n)}}{\delta}, \frac{\mathbf{X}_{(2, n)}}{\delta}, ..., \frac{\mathbf{X}_{(n, n)}}{\delta
}\right) \,\,{\mathop{=}\limits_{}^{d}} \,\,(e^{\mathbf{E}_{(1, n)}}, e^{\mathbf{E}_{(2, n)}}, ... e^{\mathbf{E}_{(n, n)}}),$$
where $\mathbf{E}_{(1, n)} \leq \mathbf{E}_{(2, n)} \leq ... \leq \mathbf{E}_{(n, n)}$ are  order statistics of independent identically distributed (i.i.d.) r.vs. with $\mathbf{E}_1 \in Exp(\alpha)$. Then, because of the multiplicative property of the exponential distribution
$$\left(\frac{\mathbf{X}_{(1, n)}}{\delta}, \frac{\mathbf{X}_{(2, n)}}{\delta}, ..., \frac{\mathbf{X}_{(n, n)}}{\delta
}\right) \,\, {\mathop{=}\limits_{}^{d}}\,\, \left(e^{\frac{1}{\alpha}\mathbf{E}_{(1, n)}^*}, e^{\frac{1}{\alpha}\mathbf{E}_{(2, n)}^*}, ... e^{\frac{1}{\alpha}\mathbf{E}_{(n, n)}^*}\right),$$
where $\mathbf{E}_{(1, n)}^* \leq \mathbf{E}_{(2, n)}^* \leq ... \leq \mathbf{E}_{(n, n)}^*$ are order statistics of i.i.d. r.vs. with $\mathbf{E}_1^* \in Exp(1)$. See e.g. de Haan and Ferreira \cite{deHaanFerreira}. Denote the logarithm with basis $e$ by log. Because of $e > 1$, $h(x) = log(x)$, is an increasing function, thus
$$log\left(\frac{\mathbf{X}_{(j, n)}}{\mathbf{X}_{(i, n)}}\right) \,\, {\mathop{=}\limits_{}^{d}}\,\, \frac{1}{\alpha}\left(\mathbf{E}_{(j, n)}^* - \mathbf{E}_{(i, n)}^*\right) \,\, {\mathop{=}\limits_{}^{d}}\,\, \frac{1}{\alpha}\mathbf{E}_{(j-i, n-i)}^*.$$
The last equality could be seen e.g. in de Haan and Ferreira \cite{deHaanFerreira} or Arnold et al. (1992) \cite{arnold1992first}.

{\bf i)} Follows by the equality $P_{\frac{E_{(j - i, n-i)}^*}{\alpha}}(x) = \alpha P_{E_{(j - i, n-i)}^*}(\alpha \, x)$, the well known relation $- \log(\rho) \stackrel{d}{=} E_{(j - i, n-i)}^*$ and the formula for probability density function (p.d.f.) of order statistics of a sample of i.i.d. r.vs. See e.g. p. 7 Nevzorov \cite{Nevzorov}.

{\bf ii)}
The mean, and the variance of the last order statistics are very well investigated. See e.g. Nevzorov \cite{Nevzorov}, p.23. Using his results and the main properties of the expectation and the variance we obtain:
\begin{eqnarray*}
  \mathbb{E}\left[log\left(\frac{\mathbf{X}_{(j, n)}}{\mathbf{X}_{(i, n)}}\right)\right] &=& \mathbb{E}\left(\frac{1}{\alpha} E_{(j - i, n-i)}^*\right) = \frac{1}{\alpha} \mathbb{E}\left( E_{(j - i, n-i)}^*\right) = \frac{1}{\alpha}\left(H_{n-i} - H_{n-j}\right)\\
 \mathbb{D}\left[log\left(\frac{\mathbf{X}_{(j, n)}}{\mathbf{X}_{(i, n)}}\right)\right] &=& \mathbb{D}\left(\frac{1}{\alpha} E_{(j - i, n-i)}^*\right) = \frac{1}{\alpha^2} \mathbb{D}\left( E_{(j - i, n-i)}^*\right) = \frac{1}{\alpha^2}\left(H_{n-i, 2} - H_{n-j, 2}\right).
\end{eqnarray*}
 \hfill Q.A.D.

In the next corollary is useful when working with finite samples. We obtain that for any $s = 2, 3, ...$, and for fixed $k = 1, 2, ...$ the estimators $Q_{k,s}$ are unbiased for $\frac{1}{\alpha}$. The accuracy of these estimators in that case is explicitly calculated.
However these estimators are applicable also for large enough samples, because  for $k \to \infty$ they are weakly consistent and asymptotically efficient.

{\bf Corollary 1.} Assume $n = k(s+1) - 1$,  $\mathbf{X}_{(1, n)} \leq \mathbf{X}_{(2, n)} \leq ... \leq \mathbf{X}_{(n, n)}$ are order statistics of independent observations on a r.v. $\mathbf{X} \in Par(\alpha, \delta)$, $\alpha > 0$, $\delta > 0$. Then, for all $s = 2, 3, ...$, and $k \in \mathbb{N}$,
 \begin{description}
   \item[i)] Denote by $\rho$ a Beta distributed with parameters $k$, and $(s-1)k$. Then
     $$Q_{k,s}  \stackrel{d}{=} \frac{-\log(\rho)}{\alpha(H_{ks-1} - H_{k-1})}  \stackrel{d}{=} \frac{E_{((s-1)k, ks-1)}}{H_{ks-1} - H_{k-1}} \stackrel{d}{=} \frac{E_{((s-1)k, ks-1)}^*}{\alpha(H_{ks-1} - H_{k-1})} ,$$
where $E_{((s-1)k, ks - 1)}$ is the $(s-1)k$-th order statistics in a sample of $ks - 1$ independent observations on i.i.d. Exponential r.vs. with parameter $\alpha$. $E_{((s-1)k, ks-1)}^*$ is the $(s-1)s$ - th order statistic of a sample of $ks-1$ independent observations on exponentially distributed r.v. with parameter $1$. Its probability density function is
 $$f_{Q_{k,s}}(x) = \frac{\alpha (H_{ks-1} - H_{k-1}) (ks-1)!}{[(s-1)k-1]!(k-1)!}(1-e^{-\alpha (H_{ks-1} - H_{k-1}) x})^{(s-1)k-1}e^{-k \alpha (H_{ks-1} - H_{k-1}) x}, \quad x > 0.$$
  \item[ii)] $\mathbb{E} Q_{k,s} = \frac{1}{\alpha},$ and
   $\mathbb{D} Q_{k,s}  = \frac{H_{ks-1, 2} - H_{k-1, 2}}{\alpha^2(H_{ks-1} - H_{k-1})^2}.$
   \item[iii)] For all $\varepsilon > 0$,
   $$P\left[\left|Q_{k,s} - \frac{1}{\alpha}\right| > \varepsilon \right] \leq \frac{H_{ks-1, 2} - H_{k-1, 2}}{\alpha^2\varepsilon^2(H_{ks-1} - H_{k-1})^2}.$$
   \item[iv)] The estimator $Q_{k,s}$ is asymptotically efficient.  For
      $k \to \infty$,
      $$\mathbb{D} Q_{k,s} \sim  \frac{H_{ks-1, 2}  - H_{k-1, 2}}{\alpha^2} \left[log\left(\frac{ks-1}{k}\right)\right]^{-2}, \quad \lim_{k \to \infty} \mathbb{D} Q_{k,s}  = 0.$$
   \item[v)] The estimator $Q_{k,s}$ is weekly consistent. More precisely, for all $\varepsilon > 0$,
   $\lim_{k\to \infty} P\left[\left|Q_{k,s} - \frac{1}{\alpha}\right| > \varepsilon \right] = 0.$
 \end{description}
{\bf Proof:} i) and ii) follow by Theorem 1, definition of $Q_{k,s}$ and the relations
$$f_{Q_{k,s}}(x) = (H_{ks-1} - H_{k-1}) f_{log\frac{\mathbf{X}_{(ks, k(s+1)-1)}}{\mathbf{X}_{(k, k(s+1)-1)}}}[x(H_{ks-1} - H_{k-1})],  \quad  \mathbb{E} Q_{k,s} = \frac{\mathbb{E}log\frac{\mathbf{X}_{(ks, k(s+1)-1)}}{\mathbf{X}_{(k, k(s+1)-1)}}}{H_{ks-1} - H_{k-1}}, \quad  \mathbb{D} Q_{k,s} = \frac{\mathbb{D}log\frac{\mathbf{X}_{(ks, k(s+1)-1)}}{\mathbf{X}_{(k, k(s+1)-1)}}}{(H_{ks-1} - H_{k-1})^2}$$

iii) is corollary of ii) and Chebyshev's inequality.

iv) It is well known that $\lim_{n \to \infty} [H_n - log(n)] = \gamma,$ where $\gamma = -\Gamma'(1) = \psi(1) \approx 0,5772$ is the Euler$-$Mascheroni constant, $\Gamma(\alpha) = \int_0^\alpha x^{\alpha-1} e^{-x} dx$, and $\psi$ is the Digamma function. By ii) for any fixed $s = 2, 3, ...$, we have
\begin{eqnarray}
  \lim_{k \to \infty} \mathbb{D} Q_{k,s} &=&  \frac{1}{\alpha^2} \lim_{k \to \infty}  \frac{H_{ks-1, 2} - H_{k-1, 2}}{(H_{ks-1} - H_{k-1})^2} = \frac{1}{\alpha^2} \lim_{k \to \infty}  \frac{H_{ks-1, 2}  - H_{k-1, 2}}{\left\{H_{ks-1} - log(ks-1) - [H_{k-1} - log(k)] + log\left(\frac{ks-1}{k}\right)\right\}^2} \\
  &=& \frac{1}{\alpha^2} \lim_{k \to \infty}  \frac{H_{ks-1, 2}  - H_{k-1, 2}}{\left[log\left(\frac{ks-1}{k}\right)\right]^2} = \frac{1}{\alpha^2}   \frac{\lim_{k \to \infty} H_{ks-1, 2}  - \lim_{k \to \infty} H_{k-1, 2}}{[log(s)]^2} = 0
\end{eqnarray}
In the last equality we have used the well known solution of the Basel problem, and more precisely the limit $\lim_{n \to \infty} H_{n,2} = \frac{\pi^2}{6}.$

v) is a consequence of ii), iii) and iv).  \hfill Q.A.D.

In the previous proof we have seen that for any fixed $s = 2, 3, ...$, $\lim_{k \to \infty} (H_{ks-1} - H_{k-1}) = log(s)$. Therefore, although
 $Q_{k,s}^*$ are biased, they are asymptotically unbiased, asymptotically normal, weakly consistent and asymptotically efficient estimators for $\frac{1}{\alpha}$. The next conclusions follow by the relation $Q_{k,s}^* = \frac{Q_{k,s} (H_{ks-1} - H_{k-1})}{log(s)}$, and the main properties of the mean and the variance.

 {\bf Corollary 2.} Assume $n = k(s+1) - 1$,  $\mathbf{X}_{(1, n)} \leq \mathbf{X}_{(2, n)} \leq ... \leq \mathbf{X}_{(n, n)}$ are  order statistics of independent observations on a r.v. $\mathbf{X} \in Par(\alpha, \delta)$, $\alpha > 0$, $\delta > 0$.
 \begin{description}
  \item[i)] Denote by $\rho$ a Beta distributed with parameters $k$, and $(s-1)k$.  Then, for all $k \in \mathbb{N}$,
  $$Q_{k,s}^*  \stackrel{d}{=} - \frac{\log(\rho)}{\alpha log(s)}  \stackrel{d}{=} \frac{E_{((s-1)k, ks-1)}}{log(s)} \stackrel{d}{=} \frac{E_{((s-1)k, ks-1)}^*}{\alpha log(s)} ,$$
where $E_{((s-1)k, ks - 1)}$ is the $(s-1)k$-th order statistics in a sample of $ks - 1$ independent observations on i.i.d. Exponential r.vs. with parameter $\alpha$. $E_{((s-1)k, ks-1)}^*$ is the $(s-1)s$ - th order statistic of a sample of $ks-1$ independent observations on exponentially distributed r.v. with parameter $1$. Its probability density function is
 $$f_{Q_{k,s}^*}(x) = \frac{\alpha\,\, log(s) (ks-1)!}{[(s-1)k-1]!(k-1)!}(1-s^{-\alpha x})^{(s-1)k-1}s^{-\alpha k x}, \quad x > 0.$$
  \item[ii)]  For all $k \in \mathbb{N}$, $\mathbb{E} Q_{k,s}^* = \frac{H_{ks-1} - H_{k-1}}{\alpha \,\,log(s)},$ and
   $\mathbb{D} Q_{k,s}^* = \frac{H_{ks-1, 2} - H_{k-1, 2}}{\alpha^2[log(s)]^2}.$
   \item[iii)] For all $\varepsilon > 0$,  and  $k \in \mathbb{N}$,
   $$P\left[\left|Q_{k,s}^* - \frac{1}{\alpha}\right| > \varepsilon \right] \leq \frac{H_{ks-1, 2} - H_{k-1, 2}}{\alpha^2\varepsilon^2[log(s)]^2}.$$
   \item[iv)]$Q_{k,s}^*$ estimator is asymptotically unbiased and asymptotically efficient.  More precisely
   $$\lim_{k\to \infty} \mathbb{E}Q_{k,s}^* = \frac{1}{\alpha}, \quad \lim_{k\to \infty} \mathbb{D} Q_{k,s}^* = 0.$$
   \item[v)] $Q_{k,s}^*$ estimator is weekly consistent. For all $\varepsilon > 0$,
   $\lim_{k\to \infty} P\left[\left|Q_{k,s}^* - \frac{1}{\alpha}\right| > \varepsilon \right] = 0.$
 \end{description}

Applications of the previous results require knowledge about confidence intervals. Therefore, in the the next theorem, we obtain asymptotic normality of these estimators which allows us later on to construct large sample confidence intervals.

 {\bf Theorem 4.} If $\mathbf{X} \in Par(\alpha, \delta)$, $\alpha > 0$, $\delta > 0$, then for all $s = 2, 3, ...$, and $k \to \infty$,
\begin{equation}\label{ANi}
 \sqrt{k(s+1) - 1}\left[log\left(\frac{X_{(ks, k(s+1) - 1)}}{X_{(k, k(s+1) - 1)}}\right)-\frac{1}{\alpha}log(s)\right] \stackrel{d}{\to} \eta_1, \quad \eta_1 \in N\left(0, \frac{s^2-1}{s\alpha^2}\right),
\end{equation}
\begin{equation}\label{ANii}
 \sqrt{k(s+1) - 1}(H_{ks-1} - H_{k-1})\left[\alpha Q_{k,s} - \frac{log(s)}{H_{ks-1} - H_{k-1}}\right] \stackrel{d}{\to} \eta_2, \quad \eta_2 \in N\left(0, \frac{s^2-1}{s}\right),
 \end{equation}
\begin{equation}\label{ANiii}\sqrt{k(s+1) - 1}\left[\alpha Q_{k,s}^* - 1\right] \stackrel{d}{\to} \eta_3, \quad \eta_3 \in N\left(0, \frac{s^2-1}{s[log(s)]^2}\right).
\end{equation}

{\bf Proof:} In this case $F^\leftarrow(p) = \frac{\delta}{\sqrt[\alpha]{1 - p}}$, and $f(x) = \frac{\alpha \delta^\alpha}{x^{\alpha + 1}}$.  Therefore,  $F^\leftarrow\left(\frac{1}{s+1}\right) = \delta \sqrt[\alpha]{\frac{s+1}{s}}$, $F^\leftarrow\left(\frac{s}{s+1}\right) = \delta \sqrt[\alpha]{s+1}$,
    $$f\left[F^\leftarrow\left(\frac{1}{s+1}\right)\right] = \frac{\alpha s^{1/\alpha + 1}}{\delta(s + 1)^{1/\alpha + 1}} \in (0, \infty),\quad f\left[F^\leftarrow\left(\frac{s}{s+1}\right)\right] = \frac{\alpha}{\delta(s + 1)^{1/\alpha + 1}}\in (0, \infty),$$

For $k \to \infty$ we have $\frac{s}{(s+1)k-1} \to \frac{1}{s+1}$, $\frac{sk}{(s+1)k-1} \to \frac{s}{s+1}$,  $\sqrt{(s+1)k-1}\left(\frac{k}{(s+1)k-1} - \frac{s}{s+1}\right) \to 0$, and $\sqrt{(s+1)k-1}$ $\left(\frac{sk}{(s+1)k-1} - \frac{s}{s+1}\right) \to 0$  therefore we can apply Smirnoff's theorem about the joint asymptotic normality of the order statistics and Theorem 1.  In order to determine $a_{F,s}$ and $b_{F,s}$ let us note that $log[F^\leftarrow(p)]' = \frac{1}{\alpha(1-p)}$. Therefore
$$a_{F,s} = \frac{1}{\{log[ F^\leftarrow(p)]\}'|_{p = \frac{1}{s+1}}} = \frac{\alpha s}{s+1}, \quad b_{F,s} = \frac{1}{\{log[ F^\leftarrow(p)]\}'|_{p = \frac{s}{s+1}}}= \frac{\alpha}{s+1}$$
The equalities
$$V = \frac{1}{(s+1)^2} \left(\frac{s}{a_{F,s}^2} - \frac{2}{a_{F,s}b_{F,s}} + \frac{s}{b_{F,s}^2}\right) = \frac{1}{\alpha^2}\left(\frac{1}{s} -\frac{2}{s} + s\right) = \frac{s^2-1}{\alpha^2 s}, \quad log\left(\frac{F^\leftarrow(\frac{s}{s+1})}{F^\leftarrow(\frac{1}{s+1})}\right)
 = \frac{1}{\alpha}log(s).$$
lead us to (\ref{ANi}). When we multiply the numerator in (\ref{ANi}) by  $\alpha[H_{ks-1} - H_{k-1}]$, and the denominator by $H_{ks-1} - H_{k-1}$, and use that $D\eta_2 = D(\alpha\eta_1) = \alpha^2 D\eta_1$ we obtain (\ref{ANii}). If we multiply both sides of (\ref{ANi}) by  $\frac{\alpha}{log(s)}$, and use that $D\eta_3 = D\left(\frac{\alpha}{log(s)}\eta_1\right) = \frac{\alpha^2}{[log(s)]^2} D\eta_1$ we obtain (\ref{ANiii}).
  \hfill Q.A.D.

 Now we are ready to compute the corresponding confidence intervals. Let us chose $\alpha_0 \in (0, 1)$ and denote by $z_{1-\frac{\alpha_0}{2}}$, $1-\frac{\alpha_0}{2}$ quantile of the standard normal distribution. Using (\ref{ANiii}), and the definition of $Q_{k,s}^*$ we obtain
 $$P\left[-z_{1-\frac{\alpha_0}{2}} \leq log(s)\sqrt{\frac{s[k(s+1) - 1]}{s^2-1}}\left(\alpha Q_{k,s}^* - 1\right)
  \leq z_{1-\frac{\alpha_0}{2}}\right] \to 1-\alpha_0, \quad k \to \infty.$$
 $$P\left[\frac{1}{Q_{k,s}^*}-\frac{z_{1-\frac{\alpha_0}{2}}}{Q_{k,s}^* log(s)} \sqrt{\frac{s^2-1}{s[k(s+1) - 1]}}\leq \alpha  \leq \frac{1}{Q_{k,s}^*} +\frac{z_{1-\frac{\alpha_0}{2}}}{Q_{k,s}^* log(s)}\sqrt{\frac{s^2-1}{s[k(s+1) - 1]}}\right] \to 1-\alpha_0, \quad k \to \infty.$$

 Therefore for any fixed $s = 2, 3, ...$, the corresponding asymptotic confidence intervals for $\alpha$ when $k \to \infty$ are:
\begin{equation}\label{ciPareto}
\left[\frac{log(s)}{log\frac{\mathbf{X}_{(ks, (s+1)k-1)}}{\mathbf{X}_{(k, (s+1)k-1)}}} -\frac{z_{1-\frac{\alpha_0}{2}}}{log\frac{\mathbf{X}_{(ks, (s+1)k-1)}}{\mathbf{X}_{(k, (s+1)k-1)}}} \sqrt{\frac{s^2-1}{s[k(s+1) - 1]}}; \frac{log(s)}{log\frac{\mathbf{X}_{(ks, (s+1)k-1)}}{\mathbf{X}_{(k, (s+1)k-1)}}}  +\frac{z_{1-\frac{\alpha_0}{2}}}{log\frac{\mathbf{X}_{(ks, (s+1)k-1)}}{\mathbf{X}_{(k, (s+1)k-1)}}}\sqrt{\frac{s^2-1
}{s[k(s+1) - 1]}}\right].
\end{equation}

{\it Simulation study}

\begin{figure}
\begin{minipage}[t]{0.5\linewidth}
    \includegraphics[scale=.56]{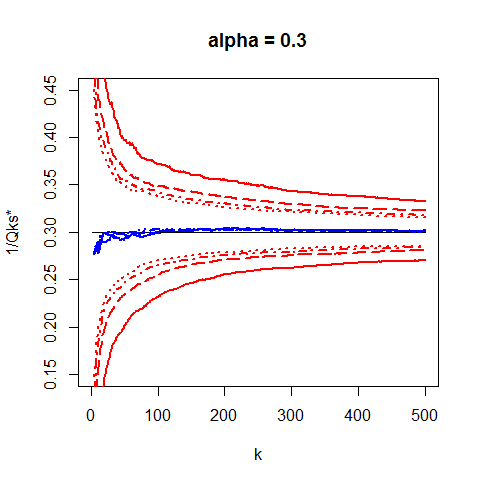}\vspace{-0.3cm}
\end{minipage}
\begin{minipage}[t]{0.49\linewidth}
    \includegraphics[scale=.56]{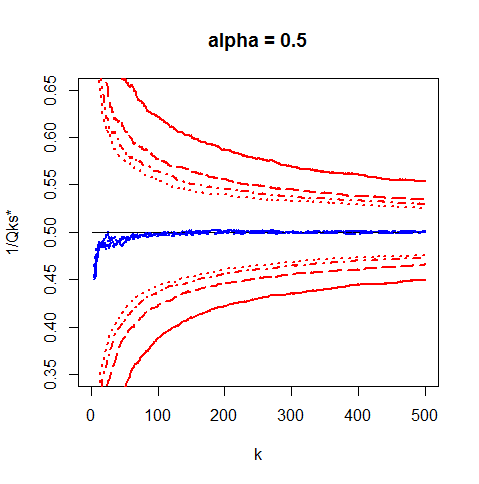}\vspace{-0.3cm}
\end{minipage}
\caption{Pareto case: Dependence of $1/Q_{k,s}^*$ and the corresponding confidence intervals(\ref{ciPareto}) on $k$, for different values of $\alpha$ and $s = 2$(solid lines), $s = 3$(dashed lines), $s = 4$(dash-dot lines), $s = 5$(dotted lines).\label{fig:Pareto0}}
\end{figure}

\begin{figure}
\begin{minipage}[t]{0.5\linewidth}
   \includegraphics[scale=.56]{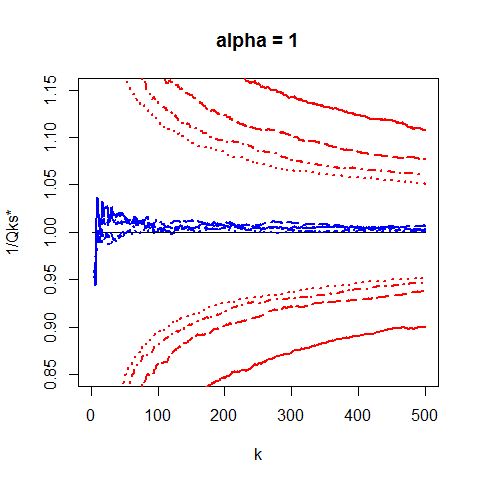}\vspace{-0.3cm}
\end{minipage}
\begin{minipage}[t]{0.49\linewidth}
    \includegraphics[scale=.56]{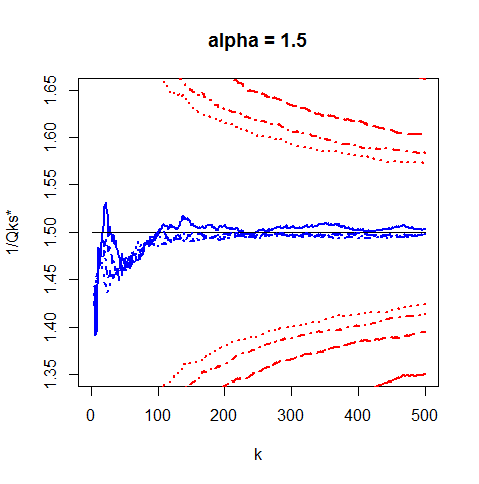}\vspace{-0.3cm}
\end{minipage}
\caption{Pareto case: Dependence of $1/Q_{k,s}^*$ and the corresponding confidence intervals(\ref{ciPareto}) on $k$, for different values of $\alpha$ and $s = 2$(solid lines), $s = 3$(dashed lines), $s = 4$(dash-dot lines), $s = 5$(dotted lines).\label{fig:Pareto1}}
\end{figure}

Let us now depict the rate of convergence of $1/Q_{k,s}^*$ for different values of $\alpha = 0.3, 0.5, 1, 1.5$ and $s = 2, 3, 4, 5$. Figures \ref{fig:Pareto0}-\ref{fig:Pareto1}  represent the dependence of $1/Q_{k,s}^*$ and the corresponding confidence intervals, on $k$. They are plotted via software R \cite{R}. The real values of $\alpha$ are plotted via straight dense line. In order to visualise the values of the estimators $\frac{1}{Q_{k,s}^*}$ for any of the lines we have simulated 100 samples of  $500(s+1)-1$ realizations of Pareto distributed r.v. with c.d.f. (\ref{Pareto}), correspondingly for $\alpha = 0.3, 0.5, 1, 1.5$. Separately for any fixed $k = 1, 2, ...$ and $s$ the values of $\frac{1}{Q_{k,s}^*}$ are averaged over these 100 samples and  presented correspondingly by dense $(s = 2)$, dashed $(s = 3)$, dash-dot $(s = 4)$, and dotted $(s = 5)$ lines. Then, for any fixed $k$, and $s$ we have computed and plotted also $0.95$-confidence intervals (red lines) for $\alpha$, calculated by formula (\ref{ciPareto}) using the averaged values of $1/Q_{k,s}^*$ instead of separate estimators $1/Q_{k,s}^*$. We observe that when $\alpha$ decreases and $s$ increases, the accuracy of the estimators improves. However, because of the sample size $k(s+1)-1$ increases with $s$ we can not chose too big $s$ for small samples.

If we compare these results with those about the well known  Hill \cite{Hill}, t-Hill\cite{Stehlik2010,JordanovaPancheva}, Deckers-Einmahl-de Haan \cite{Dekkers1989,EinmahlGuillou}, or Pickands \cite{Pickands} estimators, described very well in Embrechts et al. \cite{embrechts2013modelling}, we observe that in this case $1/Q_{k,s}^*$ estimator have better properties, especially given a small sample.

\section{CONCLUSIONS}

The paper points out good properties of couples of central order statistics  for obtaining consistent and asymptotically normal estimators of the parameter of regular variations of the tail of the c.d.f. of the observed r.v. We consider more thoroughly Pareto case, where we transform the logarithm of the fraction of the order statistics in such a way that to obtain at least asymptotically unbiased and asymptotically efficient estimator. However our results about the general case show that an analogous approach could be applied in many other cases of distributions with regularly varying tails of the c.d.f. For example: Fr\'{e}chet, Pareto, Log-logistic, Hill-horror among others. The biggest advantage of the proposed estimators is that they can be very useful for working with relatively small samples.

\section{ACKNOWLEDGMENTS}
The authors are grateful to the bilateral projects Bulgaria - Austria, 2016-2019, Feasible statistical modelling for extremes in ecology and finance, BNSF, Contract number 01/8, 23/08/2017.


\begin{thebibliography}{}

\bibitem{arnold1992first} Arnold, B.C., Balakrishnan, N., Nagaraja, H.N.: A first course in order statistics. 54, SIAM (1992)

\bibitem{Arnold2015} Arnold, B.C.: Pareto distributions. Second Edition, Chapman and Hall$/$CRC Press Taylor \& Francis Group, Boca Raton, London, New York (2015)

\bibitem{dH70} de Haan, L.: On Regular Variation and Its Application to the Weak Convergence of Sample Extremes, Mathematical Centre Tract, 32, Mathematics Centre, Amsterdam, Holland (1970)

\bibitem{deHaanStadtmueller} de Haan, L., Stadtm\"uller, U.: Generalized regular variation of second order, Journal of the Australian Mathematical Society, 61(3),381--395 (1996)

\bibitem{deHaanFerreira} de Haan, L.  and  Ferreira, A.: Extreme Value Theory: An introduction, Springer Series in Operations Research and Financial Engineering, Springer, New York (2006)

\bibitem{Dekkers1989} Dekkers, Arnold LM, Einmahl, John HJ, de Haan, Laurens: A moment estimator for the index of an extreme-value distribution, The Annals of Statistics, 1833--1855 JSTOR (1989)

\bibitem{EinmahlGuillou} Einmahl, J.H.J., Fils-Villetard, A., Guillou, A.: Statistics of extremes under random censoring, Bernoulli, 14(1), 207--227 Bernoulli Society for Mathematical Statistics and Probability (2008)

\bibitem{embrechts2013modelling} Embrechts, P., Kl{\"u}ppelberg, Cl., Mikosch, Th.: Modelling extremal events: for insurance and finance, Springer Science \& Business Media, 33,  (2013)

\bibitem{Fabian} Fabi\'an, Zd., Stehl\'\i k, M.: On robust and distribution sensitive Hill like method, IFAS res.
report.43 (2009)

\bibitem{Hill} Hill, Bruce M and others: A simple general approach to inference about the tail of a distribution, The Annals of Statistics, 5(3), 1163--1174 Institute of Mathematical Statistics(1975)

\bibitem{JordanovaPancheva} Jordanova, P.K., Pancheva, E.I.: Weak asymptotic results for t-hill estimator, Comptes rendus de l’acad{\'e}mie bulgare des sciences, 65(12), 1649--1656 (2012)

\bibitem{Jordanova2019} Jordanova, P.K.: Tails and probabilities for $p$-outliers, Submitted (2019) url={https://arxiv.org/pdf/1902.03810.pdf}

\bibitem{jordanova2017measuring} Jordanova, P.K., Petkova, M.P.: Measuring heavy-tailedness of distributions, AIP Conference Proceedings, 1910(1), 060002 AIP Publishing (2017)

\bibitem{JordanovaMilan2012} Jordanova, P., Stehl{\'\i}k, M., Fabi{\'a}n, Zd., Strelec, L.: On Estimation And Testing For Pareto Tails, Pliska Studia Mathematica Bulgarica, 22(1), 89--108 IMI-BAS, Bulgaria (2013)

\bibitem{Mandelbrot1960} Mandelbrot, B.: The Pareto-Levy law and the distribution of income, International Economic Review, 2(1), 79--106 University of Pennsylvania (1960)

\bibitem{Mosteller1946} Mosteller, F.: Qn some useful \"inefficient\" statistics, Ann. Math. Statist., 17(4), 377-408 The Institute of Mathematical Statistics (1946)

\bibitem{Nair2013} Nair, N.U., Sankaran, P.G., Balakrishnan, N.: Quantile-based reliability analysis, Birkh$\ddot{a}$user, Basel, (2013)

\bibitem{Nevzorov} Nevzorov, V.B.: Records: mathematical theory, American Mathematical Society, United States (2001)

\bibitem{Pickands} Pickands III, J.: Statistical inference using extreme order statistics, The Annals of Statistics,  119--131 JSTOR (1975)

\bibitem{Resnick87} Resnick, S.I.: Extreme Values, Regular Variation and Point Processes, Springer Series in Operations Research and Financial Engineering, Springer-Verlag, New York (1987)

\bibitem{sgouropoulos2015matching} Sgouropoulos, N., Yao, Q., Yastremiz, Cl.: Matching a distribution by matching quantiles estimation, Journal of the American Statistical Association, 110(510), 742--759 Taylor \& Francis (2015)

\bibitem{smirnov1949limit} Smirnov, N.V.: Limit distributions for the terms of a variational series, Trudy Matematicheskogo Instituta Imeni VA Steklova, 25, 3--60 Russian Academy of Sciences, Steklov Mathematical Institute of Russian Academy of Sciences (1949)

\bibitem{MultivariateDeltaMethod} Sobel, M.E.: Asymptotic confidence intervals for indirect effects in structural equation models, Sociological methodology, 13, 290-312 (1982)

\bibitem{Stehlik2010} Stehl{\'\i}k, M., Potock{\`y}, R., Waldl, H., Fabi{\'a}n, Zd.: On the favorable estimation for fitting heavy tailed data, Computational Statistics, 25(3), 485--503 Springer (2010)

\bibitem{Stehlik2012} Stehl\'\i k, M., Fabi\'{a}n, Z., St\v relec, L.: Small sample robust testing for Normality against Pareto tails, Communications in Statistics - Simulation and Computation, 41(7), 1167-1194, (2012)

\bibitem{R} R: The R Project for Statistical Computing, https://www.r-project.org/  (2019-03-28 2019)

\bibitem{Wilks1948} Wilks, S.S.: Order statistics, Bull. Amer. Math. Soc., 54(1), Part. 1, 6--50 Russian Academy of Sciences, Steklov Mathematical Institute of Russian Academy of Sciences (1948)

\end{thebibliography}
\end{document}